\documentclass[12pt]{article}
\usepackage{hyperref}
\usepackage{lscape}
\usepackage{amsmath}
\usepackage{amssymb}
\usepackage{amsthm}

\newtheorem{Theorem}{Theorem}

\newtheorem{Conjecture}{Conjecture}

\begin{document}

\begin{center}

{\Large\bf Primes in Geometric-Arithmetic Progression} \\

 
Sameen Ahmed Khan\footnote{
\url{rohelakhan@yahoo.com}, \\
Department of Engineering, 
Salalah College of Technology, 
Post Box No. 608, Postal Code 211, 
Salalah, Sultanate of Oman.
\url{http://SameenAhmedKhan.webs.com/}, \\
\url{http://sites.google.com/site/rohelakhan/}, 
\url{http://rohelakhan.webs.com/}.
}

\end{center}

\pagestyle{myheadings}
\thispagestyle{plain}
\markboth{Sameen Ahmed Khan}{Primes in Geometric-Arithmetic Progression}

\begin{abstract}
\noindent
A geometric-arithmetic progression of primes is a set of $k$ primes 
(denoted by GAP-$k$) of the form $p_1 r^j + j d$ for fixed $p_1$, $r$ 
and $d$ and consecutive $j$, 
{\it i.e}, $\left\{ p_1, \, p_1 r + d, \, p_1 r^2 + 2 d, \, p_1 r^3 + 3 d, \,  
\dots \right\}$.  We study the conditions under which, for $k \ge 2$, 
a GAP-$k$ is a set of $k$ primes in geometric-arithmetic progression.  
Computational data (along with the MATHEMATICA codes) containing progressions 
up to GAP-13 is presented. 
Integer sequences for the sets of differences $d$ corresponding to the GAPs of 
orders up to 11 are also presented. 
\end{abstract}
 
\medskip

\noindent 
{\bf Subj-class:} Number Theory \\ 

\noindent 
{\bf Mathematics Subject Classification:} 11N13 Primes in progressions \\ 

\noindent 
{\bf PACS:} Number theory, 02.10.De \\

\noindent {\bf Key words and phrases:} 
Primes, primes in arithmetic progression, 
primes in geometric-arithmetic progression, 
integer sequences. \\
 
\noindent 
{\bf Integer Sequences:} 
A000040, A000668, A040976, A172367,
A206037, A206038, A206039, A206040, A206041, A206042, A206043, A206044, A206045, 
A209202, A209203, A209204, A209205, A209206, A209207, A209208, A209209, A209210. 
 

\newpage

\section{Introduction}
Primes in arithmetic progression (denoted by AP-$k$, $k \ge 3$) refers 
to $k$ prime numbers that are consecutive terms of an 
arithmetic progression.  For example, 5, 11, 17, 23, 29 is an AP-5, a 
five-term arithmetic progression of primes with the common difference 6.  
In this example of five primes in arithmetic progression, the primes are 
not {\em consecutive} primes.  CPAP-$k$ denotes $k$ consecutive primes 
in arithmetic progression.  An example of CPAP-3 is 47, 53, 59 with the 
common difference 6.
Primes in arithmetic progression have been extensively studied both 
analytically (see the comprehensive account in~\cite{green-tao})
and numerically (see, \cite{wolfram, wikipedia}).  The largest known 
sequences contain up to 26 terms, {\it i.e}, AP-26 
and 10 consecutive primes {\it i.e}, CPAP-10 
(see ~\cite{prime-pages, ap-records} for the AP-$k$ records 
and ~\cite{prime-pages, cpap-records} for the CPAP-$k$ records).
 
The {\it geometric-arithmetic progression} refers to 
\begin{equation}\label{ga}
a, \, a r + d, \, a r^2 + 2 d, \, a r^3 + 3 d, \, \dots\,.
\end{equation}
The sequence in~(\ref{ga}) is not be confused with the 
{\it arithmetic-geometric progression}, $a$, $(a + d) r$, 
$(a + 2 d) r^2$, $(a + 3 d) r^3$, $\dots$, whose terms are composite 
by construction.  Primes in geometric-arithmetic progression is a set 
of $k$ primes (denoted by GAP-$k$) that are the consecutive terms of a 
geometric-arithmetic progression in~(\ref{ga}).  
For example 3, 17, 79 is a 3-term geometric-arithmetic progression 
({\it i.e}, a GAP-3)  
with $a = p_1 = 3$, $r = 5$ and $d = 2$.  An example of GAP-5 is, 7, 
47, 199, 911, 4423, with $p_1 = 7$, $r = 5$ and $d = 12$. 
The first term of the GAP-$k$ is called the {\em start}, 
$r$ the {\em ratio} and $d$ the {\em difference}.   
The special case of GAP-2 shall be discussed separately. 
 
For $r = 1$, GAP-$k$ reduces to AP-$k$; in this sense, GAP-$k$ is a 
generalization of the AP-$k$.  It is possible to generate GAPs with 
$p_1 = 1$, in which a case the first term of the sequence is 1 and 
has to be excluded when computing the order ($k$) as 1 is excluded 
from the set of primes.  Example of one such GAP-5 with $p_1 = 1$, 
$r = 7$ and $d = 720$ is 1, 727, 1489, 2503, 5281, 20407.  
One can also have GAPs with $p_1 = r$; an example for a GAP-5, 
is 5, 139, 353, 967, 3581, with $p_1 = r = 5$ and $d = 114$. 
There can be GAPs with composite $r$; an example of such a GAP-3 is 
7, 107, 1579 with $p_1 = 7$, a composite $r = 15$ and $d = 2$; 
and an example for GAP-5 is 11, 919, 14543, 473227, 16509011 
with $p_1 = 11$, a composite $r = 35$ and $d = 534$.  
Relevant examples are presented in Table-1 and Table-2 respectively.

\section{Results and Analysis}
The following theorem summarizes the conditions on a 
geometric-arithmetic progression to be a candidate for GAP-$k$.
 
\begin{Theorem}\label{theorem-basic}
Let GAP-$k$ denote the set of $k$ primes forming the sequence 
$\left\{p_1 r^j + j d \right\}_{j = 0}^{k-1}$, for fixed
$p_1$, $r$ and $d$.  Then it is necessary that $d$ is even; $p_1$ 
is an odd-prime coprime to $d$; $r$ is an odd-number coprime to $d$. 
When $p_1 \ne 1$ and $r \ne 1$, the maximum possible order-$k$ 
of the set is lesser of the two fixed numbers $p_1$ and the smallest 
prime factor of $r$.  
When $r = 1$, the maximum order of the set is $p_1$. 
When $p_1 = 1$, the maximum order of the set is less than the smallest 
prime factor of $r$.  
\end{Theorem}
When $d$ is odd, the alternate terms of the sequence 
$\left\{p_1 r^j + j d \right\}_{j = 0}^{k-1}$, take 
even values.  Hence, $d$ can not be odd.  When $p_1 r^j$ is even 
then again the alternate terms of the 
sequence are even. So, it is necessary that $d$ is even and  $p_1 r^j$ 
is odd, ensuring that all the terms of the sequence are odd, a 
prerequisite for them to be prime.  The first term of the sequence is 
$p_1$.  So, $p_1$ is necessarily an odd-prime.  Since, $p_1 r^j$ is odd 
it is necessary that $r$ is also odd.  For $p_1 r^j + j d$ to be prime it 
is necessary that $p_1$ and $r$ are both coprime to $d$. This proves 
the theorem except for the part related to the order of the set. 

First we consider the scenario $p_1 \ne 1$ and $r \ne 1$.  
The $(p_1 + 1)^{\text{th}}$ term of the sequence (obtained for $j = p_1$) 
is $p_1 r^{p_1} + p_1 d$, which is composite.  Hence, $k \leq p_1$. 
Let $r_1$ be the smallest prime factor of $r$.  The $(r_1 + 1)^{\text{th}}$
term of the sequence (obtained for $j = r_1$) is 
$p_1 r^{r_1} + r_1 d = p_1 r_1^{r_1} (r/r_1)^{r_1} + r_1 d$, which is composite.
Hence $k \leq r_1$.  
When $r = 1$ and $p_1 \ne 1$, the sequence simplifies to $\left\{p_1 + j d \right\}$, 
whose $(p_1 + 1)^{\text{th}}$ term is $p_1 + p_1 d$, which is composite.  
Hence, $k \leq p_1$.
When $p_1 = 1$ and $r \ne 1$, the sequence becomes $r^j + j d$, whose 
first term is 1 (for $j = 0$) and the $(r_1 + 1)^{\text{th}}$ 
term is $r^{r_1} + r_1 d = r_1^{r_1} (r/r_1)^{r_1} + r_1 d$, which is composite.  
Since, number 1 is not among the primes, $k \leq (r_1 - 1)$.  
The case $p_1 = 1$ and $r = 1$ is trivial and generates only one GAP 
(uniquely fixed with $d = 2$), which is the GAP-3, 3, 5, 7.
This completes the proof of the theorem.

For every integer, $n \ge 2$, there exists a prime $p$ such that $n < p < 2n$ 
(see for instance,~\cite{hardy}).  The elements of GAPs ($r \ne 1$) grow faster than $2n$.
Consequently, GAPs can not have consecutive primes as its members.  
Hence, we do not have consecutive primes in geometric-arithmetic progression.  

Theorem~\ref{theorem-basic} tells us the necessary conditions 
on $p_1$, $r$ and $d$, for a geometric-arithmetic progression to be a candidate 
for GAP-$k$.  
The theorem is {\it nonconstructive}, giving no clues for a recipe to generate the GAP-$k$.  
A recipe is required to choose `good' triplets $(p_1, \, r, \, d)$ 
in order to generate GAP-$k$ with larger values of $k$, and GAPs of a given order with
large number of digits. 

GAP-2 is a pair of primes of the form $(p_1, \, p_1 r + d)$ and 
consequently structurally much simpler than the larger GAP-$k$.  
For GAP-2, theorem~\ref{theorem-basic} 
simplifies to the condition that, $p_1 r$ and $d$ are coprime. 
For example, with $p_1 = 2$, $r = 6$ and $d = 5$, we have $(2, \, 17)$; 
with $p_1 = 3$, $r = 2$ and $d = 1$, $(3, \, 7)$; and 
$p_1 = 7$, $r = 100$ and $d = 211$, $(7, \, 911)$ respectively.  
In the world of primes, {\em titanic} is $1000+$ digits~\cite{prime-pages}.  
Example of a titanic GAP-2 is obtained with $p_1 = M_{4253}$, $r = 7$ 
and $d = 1422$ as $(M_{4253}, \, 7 M_{4253} + 1422)$.
Here, $M_{4253} = 2^{4253} - 1$ is the 19th Mersenne prime containing 1281 digits.  
Mersenne primes were chosen, as they are well-known and easy to express~\cite{mersenne, mersenne-oeis}. 
Pairs of primes with specific properties have been extensively studied.
For instance, {\em Sophie Germain primes} have the form $(p, \, 2p + 1)$.  
With $r = 1$ the GAP-2 further simplifies to the pair $(p, \, p + d)$.  Prime pairs 
such as {\em twin primes}, $(p, \, p + 2)$; {\em cousin primes}, $(p, \, p + 4)$;
{\em sexy primes}, $(p, \, p + 6)$, among others have been extensively studied~\cite{prime-pages}. 
 
A GAP-$k$ is said to be {\em minimal} if the minimal start $p_1$ and 
the minimal ratio $r$ are equal, {\it i.e}, $p_1 = r = p$, where $p$ is the smallest prime $\ge k$.  
Such GAPs have the form $\left\{p*p^j + j d \right\}_{j = 0}^{k-1}$. 
Minimal GAPs with different differences, $d$ do exist.  For example, the minimal GAP-5 
($p_1 = r = 5$) has the 
possible differences, 84, 114, 138, 168, $\ldots$ and the minimal 
GAP-6 ($p_1 = r = 7$) has the possible differences, 144, 1494, 1740, 2040, $\ldots$.  
A minimal GAP-$k$ is further said to be {\em absolutely minimal} if the difference $d$ is minimum.  
All the GAPs up to $k = 11$ in Table-1 are absolutely minimal.  Computationally obtained lower bounds 
of $d$ in search for higher-order minimal GAPs are also presented in Table-1.  
From theorem~\ref{theorem-basic}, it is evident that the order, $k$ of any GAP-$k$ does not 
exceed both the starting prime $p_1$ and the smallest prime factor of the ratio $r$.  
Equipped with this fact and the numerical data, we have the following two conjectures 
 
\begin{Conjecture}[Minimal Start]{}\label{conjecture-one}
The minimal starting prime, $p_1$ in a GAP-$k$ is the smallest prime $\ge k$. 
\end{Conjecture}

\begin{Conjecture}[Minimal Start and Minimal Ratio]{}\label{conjecture-two}
The minimal starting prime, $p_1$ and minimal ratio, $r$ in a GAP-$k$ is the 
smallest prime $\ge k$ and $p_1 = r$. 
\end{Conjecture}

Computational data in Table-1 supports these conjectures up to $k = 11$.  
In the context of the absolutely minimal GAPs, it is interesting to note that 
the absolutely minimal GAP-9 and the absolutely minimal GAP-10 occur for the same 
value of $d = 903030 = 31*971*(5\#) = 30101*(5\#)$, where  
$n \#$ is the primorial, $2.3.5. ... p, p \le n$. For example, $10\# = 2.3.5.7 = 210$. 
Consequently, GAP-9 is a complete 
subset of GAP-10 (in this particular instance, since they have the same $d$).  
An individual GAP-9 occurs for a higher $d = 1004250 = (5^2)*13*103*(5\#) = 33475*(5\#)$.  
This is analogous to the situation of AP-4 and AP-5 with the minimal start (which is 5).  
The corresponding sequence is $\{5 + j d \}$.  For $d = 6$ (which is the minimum difference), 
the AP-4 and AP-5 are 5, 11, 17, 23, and 5, 11, 17, 23, 29 respectively.  
The next AP-4 and AP-5 again occur at $d = 12$.  The individual AP-4 occurs only 
at $d = 18$, which is 5, 23, 41, 59. 
 
A given pair of start $p_1$ and ratio $r$, in general generates a GAP-$k$ of a certain 
order $k$ for different values of the difference $d$.  In this note, we shall focus on 
the set of differences corresponding to the minimal GAPs.
The minimal GAP-2, $\left\{2*2^j + j d \right\}_{j = 0}^{1}$ is a pair of 
primes, $(2, \, 4 + d) \equiv  (2, \, p - 4)$, where $p$ is any prime.  
Consequently, $d$ belongs to the sequence $\{p - 4 \}$, where $\{p \}$ is the 
infinite sequence of primes.  Since, the sequence $\{p \}$ is infinite, the 
sequence $\{p - 4 \}$ is also infinite.  
We shall cite various integer sequences from {\it The On-Line Encyclopedia
of Integer Sequences} (OEIS) created and  maintained by Neil Sloane.
For example, the sequence of primes, $\{p \}$ is identified by A000040 in~\cite{Sloane}.
The infinite sequence $\{p - 4 \}$: 1, 3, 7, 9, 13, 15, 19, 25, $\ldots$, is A172367~\cite{GAP-2-d}.
The integer sequences of the differences $d$, corresponding to the minimal GAPs of each 
order 3 to 11 are presented in~\cite{GAP-3-d}-\cite{GAP-11-d}.  
In general, there are no reasons to believe that the sequence of the differences $d$ 
corresponding to any GAP (minimal or non-minimal) are finite. 
Analogous sequences for the differences also exist for the primes in arithmetic progression. 
See~\cite{AP-2-d}-\cite{AP-differences} for the sequences of differences 
corresponding to the primes in arithmetic progression with the minimal start.
A study of these integer sequences may provide a pattern, which will potentially guide us in our
search for higher order GAPs and APs. 
 
From the computed sequences, we note that the set of differences for a given 
minimal GAP-$k$ have a {\em common} $k$-dependent multiplicative factor.  
This factor has been indicated as $(...\#)$ in Table-1.  We have included only the 
common factor. Individual differences $d$ do have additional factors.  For instance, 
the first difference for the minimal GAP-11 is $443687580 = 2112798 (7\#) = 14789586(5\#)$. 
The second difference is not a multiple of $(7\#)$ and hence we have shown the first 
difference as $14789586(5\#)$ in Table-1. 
Theorem~\ref{theorem-basic} restricts the values of the differences $d$ for 
any GAP-$k$ ($k \ge 3$) to be even ({\it i.e} multiples of 2).
There are additional restrictions on the values of the differences $d$ for 
minimal GAPs of a given order, as shown below   

\begin{Theorem}[Factors of $d$]{}\label{theorem-primorial} ~~ 

\begin{enumerate}
{\itemsep -0.103cm
\item
The values of the differences $d$ for each of the minimal GAP-$k$, $k \ge 5$ are 
multiples of a $k$-dependent factor denoted along with the order $k$ by $(k:~...\#)$, 
where $\#$ is the primorial.  They are 
$(5:~3\#)$, 
$(6-7:~3\#)$, 
$(8-11:~5\#)$, 
$(12-13:~7\#)$, 
$(14-17:~5\#)$, 
$(18:~7\#)$, 
$(19:~11\#)$, 
$(20-23:~11\#)$, 
$(24-29:~13\#)$, 
$(30-31:~13\#)$, 
$(32-37:~19*11\#)$, 
$(38-41:~13\#)$, 
$(42:~17\#)$, 
$(43:~19\#)$, 
$(44-47:~23*17\#)$, 
$(48-53:~17\#)$, 
$(54:~29*13\#)$, 
$(55-58:~29*19*13\#)$, 
$(59:~29*19\#)$, 
$(60-61:~31*19*13\#)$, 
$(62-67:~31*19\#)$, 
$(68-71:~17\#)$, 
$(72-73:~37*23*17\#)$, 
$(74-79:~23\#)$, 
$(80-83:~41*19\#)$, 
$(84-89:~31*23\#)$, 
$(90-97:~23\#)$, 
$(98-99:~23\#)$, 
$(100-101:~31*23\#)$, 
$(102-103:~23\#)$ respectively.

\item
The values of the differences $d$ for all minimal GAP-$k$, $k \ge 5$ are multiples of $(3\#)$. 
\item
The values of the differences $d$ for all minimal GAP-$k$, $k \ge 8$ are multiples of $(5\#)$. 
\item
The values of the differences $d$ for all minimal GAP-$k$, $k \ge 18$ are multiples of $(7\#)$. 
\item
The values of the differences $d$ for all minimal GAP-$k$, $k \ge 19$ are multiples of $(11\#)$. 
\item
The values of the differences $d$ for all minimal GAP-$k$, $k \ge 38$ are multiples of $(13\#)$.
}
\end{enumerate} 
\end{Theorem}
The proof of the theorem~\ref{theorem-primorial} is based on modular arithmetic and is 
presented in Appendix-A.  Factors up to $k = 103$ are presented in Table-3. 
In an arithmetic progression (AP-$k$) with the minimal start, the pattern of the differences 
is known to be a multiple of $k\#$, where $k$ is the largest prime $\le k$ (if $k$ is not a prime).  
If $k$ is a prime than the common difference is a multiple of $(k - 1)\#$  
(see \cite{wolfram, wikipedia}). 
Unlike in the case of the AP-$k$ with the minimal start, there is no obvious pattern in the case 
of the minimal GAP-$k$.  The factor $(...\#)$ is not even monotonic.  
In none of the cases, it has been possible to establish the factors containing 
the higher powers of 2, 3, 5, or 7.  Theorem~\ref{theorem-primorial} only gives the restrictions 
on the common difference $d$ in order for the generating sequence to be a candidate for minimal GAP.
The existence of the minimal GAP-$k$, $k \ge 12$ is yet to be established (numerically or otherwise).
The extension of the theorem~\ref{theorem-primorial} to non-minimal GAPs is also discussed in 
Appendix-A.
 
So far, we have considered the GAPs from the sequence $\left\{p_1 r^j + j d \right\}_{j = 0}^{k-1}$.
The sequence, $\left\{p_1 r^j + j d \right\}$ can have sets of primes 
for consecutive $j$, not necessarily starting with $j = 0$. 
For example, the sequence, $\left\{5*5^j + 4j \right\}_{j = 7}^{j = 9}$ generates the 
GAP-3, 390653, 1953157, 9765661.  
Another example is the sequence, $\left\{13*13^j + 156497*(11\#)j \right\}_{j = 3}^{j = 12}$ 
generating the GAP-10, 
1084552771, 1446403573, 1812367159, 2231796937, 3346287211, 13496563933, 141112064479, 
1795775474737, 23302061711251, 302879444689093. 
Such sets, not starting with $j = 0$, can not be put in the form
$\left\{P_1 R^j + j D \right\}_{j = 0}^{j = k-1}$, where $P_1$, $R$ and $D$ are derived from
$p_1$, $r$ and $d$.
Importantly, the orders of such GAPs are restricted, as we shall soon see.
 
\begin{Theorem}\label{theorem-free}
Let $\left\{p_1 r^j + j d \right\}_{j = 0}^{p^{\prime} - 1}$ be a 
GAP-$p^{\prime}$ of order $p^{\prime}$, where $p^{\prime}$ is smaller of the two 
primes $p_1$ and $r_1$ the smallest prime factor of $r$.  
Then the infinite sequence $\left\{p_1 r^j + j d \right\}_{j = p^{\prime}}^{\infty}$ does not 
have any GAPs of order $\ge p^{\prime}$.  
\end{Theorem}
Since, theorem~\ref{theorem-basic} requires $p_1$ to be an odd-prime and $r$ 
to be any odd number, $p^{\prime} \ge 3$.  When $r = 1$, the GAP is reduced 
to an AP and we have $p^{\prime} = p_1 \ge 3$.  
Let us recall that theorem~\ref{theorem-basic} forbids GAPs of orders greater 
than $p^{\prime}$ throughout the interval $[j = 0, \, \infty)$.
While proving theorem~\ref{theorem-basic}, we saw that 
$(p_1 + 1)^{\text{th}}$ and the $(r_1 + 1)^{\text{th}}$ terms of the 
sequence, $\left\{p_1 r^j + j d \right\}_{j = 0}^{\infty}$ are composite. 
The $(n p_1 + 1)^{\text{th}}$ term (where, $n = 1, 2, 3, \ldots$) of the sequence 
(obtained for $j = n p_1$) is $p_1 r^{n p_1} + n p_1 d$.  This term is composite and 
belongs to the interval $[j = p_1, \, \infty)$.
There are only $(p_1 - 1)$ terms between any two successive $(n p_1 + 1)^{\text{th}}$ 
and $(\overline{n+1} p_1 + 1)^{\text{th}}$ terms.  So, the interval $[j = p_1, \, \infty)$
can not have any GAPs of order more than $(p_1 - 1)$.
The $(n r_1 + 1)^{\text{th}}$ term (obtained for $j = n r_1$) of the sequence is 
$p_1 r^{n r_1} + n r_1 d = p_1 r_1^{n r_1} (r/r_1)^{n r_1} + r_1 d$.  This is also composite.
Similar arguments forbid the GAPs of order more than $(r_1 - 1)$ in the 
interval $[j = p_1, \, \infty)$.  This proves the theorem.
In passing we note that, when $r = 1$, the GAP is reduced to an AP.  Then 
theorem~\ref{theorem-free} tells us that the 
sequence, $\left\{p_1 + j d \right\}_{j = p_1}^{\infty}$ does not have any APs of order $\ge p_1$.  
 
Theorem~\ref{theorem-free} is for the restricted case of GAPs, whose order is any 
prime $p^{\prime} \ge 3$  
and forbids the existence of GAPs of order $p^{\prime}$ in the infinite interval $[j = p^{\prime}, \, \infty)$.
Moreover, the theorem is silent about the absence (or existence) of GAPs of orders lower than 
$(p^{\prime} - 1)$ in the interval $[j = p_1, \, \infty)$.  It is interesting to note that 
the nine sequences $\left\{p*p^j + j d \right\}_{j = k}^{1000}$, (for $k = 3$ to 11), 
for the choice of the absolutely minimal triplets $(p, \, p, \, d)$ in Table-1 do not have 
any GAPs of order 3 to $p$ respectively.  Some of them do have one or two GAP-2 in the 
interval, $[j = k, \, 1000]$ respectively.    
This numerical data provides room for extending the theorem~\ref{theorem-free} to cases, when the 
order of minimal GAP is not a prime.  This leads to the conjecture 

\begin{Conjecture}[GAP free]{}\label{conjecture-three}
Whenever the sequence $\left\{p*p^j + j d \right\}_{j = 0}^{k - 1}$ has the minimal GAP-$k$,
the rest of the infinite sequence, $\left\{p*p^j + j d \right\}_{j = k}^{\infty}$ 
does not have any GAPs of orders $\ge 3$. 
\end{Conjecture}

\section{Concluding Remarks}
For a given triplet, $(p_1, \, r, \, d)$, the 
sequence, $\left\{p_1 r^j + j d \right\}_{j = 0}^{N}$,  
may not always generate very many primes. 
For example, the sequence, $\left\{5*3^j + 2j \right\}$, takes prime values for,  
$j = 0$, 1, 7, 29, 49, 83, 436, 536, 1274, $\dots$.  
The sequence,  $\left\{7*13^j + 36j \right\}_{j = 0}^{1000}$, has only a 
single pair ({\it i.e}, a GAP-$2$) for $j = 0, 1$, which is $(7, \, 127)$.
Similar is the situation for 
a wide range of $(p_1, \, r, \, d)$, making it very hard to find GAPs. 
GAP-3 and GAP-4 with 159 digits were obtained using the Mersenne primes. 

Numerical data in this article was computed initially 
(up to GAP-6 in Table-1), using the Microsoft {\it EXCEL}~\cite{EXCEL}.  
The primality of the numbers generated by EXCEL was checked using the database of 
primes at {\it The Prime Pages}~\cite{prime-pages} and the {\em Sequence A000040}, from 
{\it The On-Line Encyclopedia of Integer Sequences} (OEIS), created and 
maintained by Neil Sloane~\cite{Sloane}.  
For higher orders, we are using the {\it MATHEMATICA}~\cite{MATHEMATICA}.
Search for GAPs with ever larger $k$ and  
geometric-arithmetic progressions containing larger primes is in progress.
 
We end this note with several open questions, similar to the ones, which exist for the 
primes in arithmetic progression~\cite{green-tao}.  
Are there arbitrarily long geometric-arithmetic progressions of primes?  
Are there infinitely many $k$-term geometric-arithmetic progressions consisting of $k$ primes?  
Do the prime numbers contain infinitely many geometric-arithmetic progressions of length $k$ for all $k$? 
Are there are infinitely many GAP-$k$ for any $k$?  
We conjecture that the answer to all the above questions is in the affirmative. 
 
\newpage

\begin{center}
{\Large\bf 
Appendix-A: \\ Proof of Theorem~\ref{theorem-primorial}}
\end{center}
 
Theorem~\ref{theorem-basic} states that the common factor of all differences of 
any GAP-$k$, $k\ge 3$ is 2 ({\it i.e}, $d$ is even).  
Theorem~\ref{theorem-primorial} states additional factors of all the differences 
$d$ of a given minimal GAP-$k$, $k \ge 5$.  The theorem essentially consists of two parts:
part-1 is for the specific GAP-$k$ with $k$ up to 103; part-2 has global statements giving 
the common factor of all differences of any minimal GAP, whose order exceeds a particular 
number.  Its proof is based on modular arithmetic.  As we shall soon see, it suffices to 
demonstrate the procedure of proving the statements in a few specific cases of both part-1
and part-2 respectively.  Rest of the statements in the theorem for higher orders can be 
proved closely following the procedures established for lower orders.  In fact the 
methodology presented can be used to derive results and extend the theorem to still higher 
orders and importantly to the non-minimal GAPs. 

The minimal GAP-5 is defined by the sequence $\left\{5*5^j + j d \right\}_{j = 0}^{4}$ and 
the common difference $d$ is restricted in such a way that the defining sequence has 5 primes. 
The first four terms of this sequence belong to GAP-4.   
The residues of this sequence $\pmod {3}$ are $\{2, \, 1 + d, \, 2 + 2 d, \, 1, \, 2 + d \}$. 
Primality requires the second residue $1 + d$ such that $d \not\equiv 2 \pmod {3}$ and the 
fifth residue $2 + d$ such that $d \not\equiv 1 \pmod {3}$.  Consequently, $d \equiv 0 \pmod {3}$.  
Otherwise, the second and fifth terms in the defining sequence would be multiples of 3 and not prime.
Hence, $d$ is {\em necessarily} a multiple of 3.  
Theorem~\ref{theorem-basic} restricts the values of the differences $d$ of any GAP-$k$, $k\ge 3$ to 
be multiples of 2.  Consequently, the values of the differences for GAP-5 are restricted to be multiples
of $(3\#)$.  
The fifth residue does not belong to GAP-4 and hence the result is not applicable to GAP-4.
The third residue $2 + 2 d$ is degenerate as it gives the same information as the second residue.
Among the five residues, the first was {\em numeric} ({\it i.e}, free of $d$) and the third was degenerate.  

The defining sequence for the minimal GAP-7 is $\left\{7*7^j + j d \right\}_{j = 0}^{6}$. 
The corresponding residues $\pmod {3}$ are $\{1, \, 1 + d, \, 1 + 2 d, \, 1, \, 1 + d, \, 1 + 2 d, \, 1 \}$.
Since, GAP-6 is defined by the same sequence except for the index, its residues are the same 
as the first six residues for GAP-7.  The second and third residues are sufficient to establish
that $d$ is a multiple of 3 for both GAP-6 and GAP-7.  
Consequently, the values of the differences for GAP-6 and GAP-7 are restricted to be multiples
of $(3\#)$.  
The residues $\pmod {5}$ are $\{2, \, 4 + d, \, 3 + 2 d, \, 1 + 3 d, \, 2 + 4 d, \, 4, \, 3 + d \}$.  
The first and sixth residues are numeric. 
The second, fourth and fifth residues require $d \not\equiv 1 \pmod {5}$, $d \not\equiv 3 \pmod {5}$, 
and $d \not\equiv 2 \pmod {5}$ respectively.  The third and seventh residues are degenerate.  
The case, $4 \pmod {5}$ remains unaddressed and hence $d$ need not be multiple of 5. 
 
The presence of numeric and degenerate residues of a given defining sequence hinders the 
larger factors.  
Rest of the results in part-1 of the theorem are proved closely following the procedure 
used for GAP-4 to GAP-7.  The procedure is straightforward but becomes laborious 
as the order-$k$ grows.  We have used the MATHEMATICA to compute the residues~\cite{MATHEMATICA}.  
Following are the residues for GAP-7  $\pmod {5}$ 
{\tt
\begin{verbatim}
In[1]:= Clear[p];
p = 7;
PolynomialMod[{p, p*p + d, p*p^2 + 2*d, p*p^3 + 3*d, 
	p*p^4 + 4*d, p*p^5 + 5*d, p*p^6 + 6*d}, 5]

Out[3]= {2, 4 + d, 3 + 2 d, 1 + 3 d, 2 + 4 d, 4, 3 + d}
\end{verbatim}
}

In part-1 of the theorem, the results are for individual GAP-$k$, $k \le 103$.  We have 
demonstrated the procedure up to $k = 7$ and it is straightforward to extend it to higher 
orders.  Part-2 onwards the statements are global and the procedure is as follows. 
The differences for the minimal GAP-5 are divisible by 3.  Since, $5 \equiv 2 \pmod {3}$ 
the result is applicable to all those primes $> 5$, whose residues $\pmod {3}$ are 2.
The differences for GAP-7 are divisible by 3 and $7 \equiv 1 \pmod {3}$.  The result is 
again applicable to all those primes $> 7$, whose residues $\pmod {3}$ are 1.  The 
result was individually proved for GAP-6 in part-1.  
The non-zero residues of 3 are $\{1, 2 \}$.  Consequently the differences for 
all GAP-$k$, $k \ge 5$ are divisible by $(3\#)$ with the factor 2 coming from 
theorem~\ref{theorem-basic}.

The non-zero residues of 5 are $\{1, 2, 3, 4 \}$, and the corresponding 4 primes 
with these residues are $\{11, 17, 13, 19 \}$.  Note, that $7 \equiv 2 \pmod {5}$ 
but its differences are not divisible by 5 as seen in part-1.
We now individually examine the four generating sequences for GAP-11, GAP-17, GAP-13 
and GAP-19 respectively and conclude that the differences for each of them are multiples of 5. 
The factor of $(3\#)$ is already established, so we conclude that all GAP-$k$, $k\ge 19$ have 
their differences as multiples of $(5\#)$.  The inequality $k \ge 19$, is refined by using the 
results in part-1.  The factor of $(5\#)$ was established for the lower order 
GAP-$k$, $k = 8$ to $k = 18$ in part-1. 
Hence, the values of the differences $d$ for all minimal GAP-$k$, $k \ge 8$ are multiples of $(5\#)$.
The set $\{11, 17, 13, 19 \}$ was only a {\em candidate} set.  Had the differences for any of 
the GAP-11, GAP-17, GAP-13 or GAP-19 failed to be divisible by 5, we would have examined the 
GAPs of higher orders, corresponding to that residue.  

The non-zero residues of 7 are $\{1, 2, 3, 4, 5, 6 \}$ and the corresponding 6 primes 
are $\{29, 23, 31, 53, 19, 13 \}$.  The primes $11 \equiv 4 \pmod {7}$ and $17 \equiv 3 \pmod {7}$ 
are not relevant in view of the results in part-1.  Following the procedure used for establishing 
the factors $(3\#)$ and $(5\#)$, we conclude that the differences $d$ for all 
minimal GAP-$k$, $k \ge 18$ are multiples of $(7\#)$. 
 
The candidate set of 10 primes corresponding to the non-zero residues of 11 
is $\{23, 79, 47, 37, 71, 61, 29, 19, 31, 43 \}$.  The differences $d$ for each of the
GAPs of these orders are divisible by $(11\#)$.  
The largest prime in this set is 79 and a spontaneous result is that all GAP-$k$, $k \ge 79$ 
have their differences as multiples of $(11\#)$.  
Using the results from part-1, we refine the inequality and conclude that
the differences $d$ for all minimal GAP-$k$, $k \ge 19$ are multiples of $(11\#)$. 
 
The candidate set of 12 primes corresponding to the non-zero residues of 13 
is $\{53, 41, 29, 43, 31, 71, 59, 47, 61, 101, 89, 103 \}$.  The largest prime in 
this set is 103.  Hence, the results in part-1 are up to $k = 103$.  The candidate set 
successfully works and we conclude that the differences $d$ for all 
minimal GAP-$k$, $k \ge 38$ are multiples of $(13\#)$. 
 
The candidate set of 16 primes corresponding to the non-zero residues of 17 
is $\{103, 53, 71, 89, 73, 193, 109, 127, 43, 163, 79, 97, 47, 167, 83, 101 \}$.  
The largest prime in this set is 193.

The candidate set of 18 primes corresponding to the non-zero residues of 19 
is $\{191, 59, 79, 61, 43, 101, 83, 103, 199, 67, 163, 107, 89, 109, 167, 149, 
131, 37 \}$.  
The largest prime in this set is 199. 
  
The procedure of proving this theorem can be applied to the non-minimal GAPs.  
The common factor $...\#$ of the differences $d$ of a GAP-$k$ with the start $p_1$, 
and ratio $r$ shall be denoted by $(k: p_1, \, r, \, ...\#)$.  
The examples are 
$(3: 5, 7, 3\#)$m
$(3: 2^{521} - 1, 19, 3\#)$, 
$(4: 11, 35, 2\#)$, 
$(4: 2^{521} - 1, 5, 2\#)$, 
$(5: 47, 2^{31} - 1, 3\#)$, 
$(7: 7, 11, 5\#)$, 
$(7: 7, 13, 3\#)$,
$(7: 7, 17, 3\#)$,
$(7: 7, 19, 3\#)$,
$(7: 11, 7, 3\#)$, 
$(7: 11, 13, 3\#)$,
$(7: 11, 17, 3\#)$,
$(7: 13, 7, 3\#)$,
$(7: 17, 7, \#)$,
$(7: 19, 7, 3\#)$,
$(7: 19, 23, 3\#)$,
$(7: 23, 19, 3\#)$,
$(11: 11, 13, 7\#)$ and 
$(11: 13, 11, 5\#)$.  
The choice of the $(k: p_1, \, r, \, ...\#)$ in the above examples includes the cases 
covered in Table-2.
 
\newpage

\begin{center}
{\Large\bf 
Appendix-B: \\ MATHEMATICA Codes}
\end{center}
 
Most of the data in this article was computed using the versatile package {\it MATHEMATICA}~\cite{MATHEMATICA}. 
The following program searches for the values of the differences $d$ for the minimal GAP-5, in the 
range $[0, 1000]$.     
{\tt
\begin{verbatim}
In[1]:=	Clear[p]; p = 5;
gapset5d = {};
Do[If[PrimeQ[{p, p*p + d, p*p^2 + 2*d, p*p^3 + 3*d, 
     p*p^4 + 4*d}] == {True, True, True, True, True}, 
  AppendTo[gapset5d, d]], {d, 0, 10^3}]; gapset5d // Timing
 
Out[4]= {6.50521*10^-19, {84, 114, 138, 168, 258, 324, 348, 
	 		462, 552, 588, 684, 714, 744, 798, 882, 894, 972}}
\end{verbatim}
}
The output is the set of 17 values of $d$: $\{84, 114, 138, \dots, 972 \}$.
The above program (christened {\em runner}) picks the values of $d$ but skips 
the finer details.  The following program (christened {\em walker}) gives the 
complete GAP sets corresponding to each $d$.   
{\tt
\begin{verbatim}
In[5]:= f[n_, m_] := (5)*(5)^n + n*m;
	Column[Table[{m, 
   Cases[Table[{f[n, m], f[n + 1, m], f[n + 2, m], f[n + 3, m], 
       	f[n + 4, m]}, {n, 0, 5}], {a1_, a2_, a3_, a4_, a5_} /; 
      	PrimeQ[{a1, a2, a3, a4, a5}] == {True, True, True, True, 
        True}]}, {m, 114, 114}]] // Timing

Out[6]= {4.33681*10^-19, {114, {{5, 139, 353, 967, 3581}}}
\end{verbatim}
}
In the above program we choose the difference 114 and obtained the corresponding GAP-5: 
$\{5, 139, 353, 967, 3581 \}$.  As the name goes the walker is much slower than the 
runner and hence, not suitable for generating the sequence of differences $d$.  
The runner can be made into an {\em accelerator} by replacing $\{d, 0, 10^{\wedge} 3 \}$ 
with $\{d, 0, 10^{\wedge} 3, \, {\mbox{\boldmath $2$}} \}$ and confining the search to 
multiples of 2 (as restricted by theorem~\ref{theorem-basic}).  It can be further 
accelerated by the replacement of $\{d, 0, 10^{\wedge} 3, \, {\mbox{\boldmath $2$}} \}$ 
with $\{d, 0, 10^{\wedge} 3, \, {\mbox{\boldmath $6$}} \}$ and refining the search to 
multiple of $(3\#) = 6$ (as restricted by theorem~\ref{theorem-primorial}).  Such 
replacements are relevant as the numbers grow.  It is straightforward to extend the 
above programs (for GAP-5) to higher orders. 

\newpage 

\begin{landscape}
 
\noindent
{\bf Table-1:} 
Primes in Geometric-Arithmetic Progression with 
{\em minimal start}, $p_1$, {\em minimal ratio}, $r$ and the {\em minimal difference}, $d$. 

\begin{center}
\begin{tabular}{|c|c|c|c|c|c|c|}
\hline
$k$ & $p_1$ & $r$ & $d$ & 
$\begin{array}{c} {\rm Primes~ of~ the~ form}, \\ 
p_1 *r^n + n d, 
{\rm for} ~$n = 0$ ~{\rm to} ~k - 1
\end{array}$ & $\begin{array}{c} {\rm Digits}\\ {\rm of~First}\end{array}$ 
& $\begin{array}{c} {\rm Digits}\\ {\rm of~Last}\end{array}$ \\
\hline
2 & 2 & 2 & 1 & $2*2^n + n$ & 1 & 1 \\
\hline	
3 & 3 & 3 & 2 & $3*3^n + 2n$ & 1 & 2 \\
\hline
4 & 5 & 5 & $3(2\#)$ & $5*5^n + 3(2\#)n$ & 1 & 3 \\
\hline
5 & 5 & 5 & $14(3\#)$ & $5*5^n + 14(3\#)n$ & 1 & 4 \\
\hline
6 & 7 & 7 & $24(3\#)$ & $7*7^n + 24(3\#)n$ & 1 & 6 \\
\hline
7 & 7 & 7 & $554(3\#)$ & $7*7^n + 554(3\#)n$ & 1 & 6 \\ 
\hline
8 & 11 & 11 & $2087(5\#)$ & $11*11^n + 2087(5\#)n$ & 2 & 9 \\
\hline
9 & 11 & 11 & $30101(5\#)$ & $11*11^n + 30101(5\#)n$ & 2 & 10 \\
\hline
10 & 11 & 11 & $30101(5\#)$ & $11*11^n + 30101(5\#)n$ & 2 & 11 \\
\hline
11 & 11 & 11 & $14789586(5\#)$ & $11*11^n + 14789586(5\#)n$ & 2 & 12 \\
\hline
$12-13$ & 13 & 13 & $> {25*10^7} \times (7\#)$ & & & \\
\hline
$14-17$ & 17 & 17 & $> {6*10^7} \times (5\#) $ & & & \\
\hline
$18$ & 19 & 19 & $> {10*10^7} \times (7\#) $ & & & \\ 
\hline
$19$ & 19 & 19 & $> {10*10^7} \times (11\#) $ & & & \\
\hline
$20-23$ & 23 & 23 & $> {5*10^7} \times (11\#)$ & & & \\
\hline
$24-29$ & 29 & 29 & $> {6*10^7} \times (13\#)$ & & & \\
\hline
$30-31$ & 31 & 31 & $> {11*10^7} \times (13\#)$ & & & \\
\hline	
\end{tabular}
\end{center}

\noindent
$n \#$ is the primorial, $2.3.5. ... p, p \le n$. For example, $10\# = 2.3.5.7 = 210$. \\
The symbol $>$ indicates the computationally obtained lower bound of 
the {\em minimal difference}, $d$ in search for the GAPs of the corresponding orders. 

\newpage
 
\noindent
{\bf Table-2:} Miscellaneous examples of Primes in Geometric-Arithmetic Progression 
 
\begin{center} 
\begin{tabular}{|c|c|c|c|c|c|c|}
\hline
$k$ & $p_1$ & $r$ & $d$ & 
$\begin{array}{c} {\rm Primes~ of~ the~ form}, \\ 
p_1 *r^n + n d, 
{\rm for} ~$n = 0$ ~{\rm to} ~k - 1
\end{array}$ & $\begin{array}{c} {\rm Digits}\\ {\rm of~First}\end{array}$ 
& $\begin{array}{c} {\rm Digits}\\ {\rm of~Last}\end{array}$ \\
\hline
2 & 2 & 5 & 7 & $2*5^n + 7$ & 1 & 2 \\
\hline	
2 & 13 & 80 & 53 & $13*80^n + 53$ & 1 & 4 \\
\hline
2 & $2^{4423} - 1$ & 7 & 802 & $(2^{4423} - 1)*7^n + 802n$ & 1332 & 1333 \\
\hline
3 & 5 & 7 & $4(3\#)$ & $5*7^n + 4(3\#)n$ & 1 & 3 \\
\hline
3 & $2^{127} - 1$ & 3 & 7390 & $(2^{127} - 1)*3^n + 7390n$ & 39 & 40 \\
\hline
3 & $2^{521} - 1$ & 3 & 1106 & $(2^{521} - 1)*3^n + 1106n$ & 157 & 158 \\
\hline
3 & $2^{521} - 1$ & 19 & $4365(3\#)(2^{127} - 1)$ & $(2^{521} - 1)*19^n + 4365(3\#)(2^{127} - 1)n$ & 157 & 160 \\
\hline
4 & 11 & 35 & $12(2\#)$ & $11*35^n + 12(2\#)n$ & 1 & 15 \\
\hline
4 & $2^{521} - 1$ & 5 & $33936(2\#)$ & $(2^{521} - 1)*5^n + 33936(2\#)n$ & 157 & 159 \\
\hline
5 & 47 & $2^{31} - 1$ & $13554(3\#)$ & $47*(2^{31} - 1)^n + 13554(3\#)n$ & 2 & 39 \\
\hline
10 & 17 & 17 & 7700(5\#) & $17*17^n + 7700(5\#)n$ & 2 & 13 \\
\hline
11 & 13 & 11 & $129262(5\#)$ & $13*11^n + 129262(5\#)n$ & 2 & 12 \\
\hline
11 & 13 & 13 & $983234(7\#)$ & $13*13^n + 983234(7\#)n$ & 2 & 13 \\
\hline
11 & 103 & 103 & 2900641(23\#) & $103*103^n + 2900641(23\#)n$ & 3 & 23 \\
\hline	
12 & 19 & 19 & $ 9345011(11\#)$ & $19*19^n + 9345011(11\#)n$ & 2 & 16 \\
\hline
12 & 31 & 31 & $70167912(13\#)$ & $31*31^n + 70167912(13\#)n$ & 2 & 18 \\
\hline	
13 & 19 & 19 & $ 9345011(11\#)$ & $19*19^n + 9345011(11\#)n$ & 2 & 17 \\
\hline	
\end{tabular}
\end{center}
 
\noindent
$M_{31} = 2^{31} - 1$, $M_{127} = 2^{127} - 1$, $M_{521} = 2^{521} - 1$ 
and $M_{4423} = 2^{4423} - 1$ are Mersenne primes~\cite{mersenne, mersenne-oeis}. 
 
\end{landscape}
 
\newpage

\noindent
{\bf Table-3:} 
The differences $d$ for the minimal GAP of each order are multiples of a common factor 
 
\begin{center} 
\begin{tabular}{|c|c|c|}
\hline
Order-$k$ & Generating Prime $p$ & Common Factor \\
\hline
2 & 2 & 1 \\
\hline
3 & 3 & 2 \\
\hline
4 & 5 & 2 \\
\hline
5 & 5 & $3\#$ \\
\hline
6--7 & 7 & $3\#$ \\
\hline
8--11 & 11 & $5\#$ \\
\hline	
12--13 & 13 & $7\#$ \\
\hline
14--17 & 17 & $5\#$ \\
\hline
18 & 19 & $7\#$ \\
\hline
19 & 19 & $11\#$ \\
\hline
20--23 & 23 & $11\#$ \\
\hline
24--29 & 29 & $13\#$ \\
\hline
30--31 & 31 & $13\#$ \\
\hline
32--37 & 37 & $19*11\#$ \\
\hline
38--41 & 41 & $13\#$ \\
\hline
42 & 43 & $17\#$ \\
\hline
43 & 43 & $19\#$ \\
\hline
44--47 & 47 & $23*17\#$ \\
\hline
48--53 & 53 & $17\#$ \\
\hline
54 & 59 & $29*13\#$ \\
\hline
55--58 & 59 & $29*19*13\#$ \\
\hline
59 & 59 & $29*19\#$ \\
\hline
60--61 & 61 & $31*19*13\#$ \\
\hline
62--67 & 67 & $31*19\#$ \\
\hline
68--71 & 71 & $17\#$ \\
\hline
72--73 & 73 & $37*23*17\#$ \\
\hline
74--79 & 79 & $23\#$ \\
\hline
80--83 & 83 & $41*19\#$ \\
\hline
84--89 & 89 & $31*23\#$ \\
\hline
90--97 & 97 & $23\#$ \\
\hline
98--99 & 101 & $23\#$ \\
\hline
100--101 & 101 & $31*23\#$ \\
\hline
102--103 & 103 & $23\#$ \\
\hline
\end{tabular}
\end{center}

\newpage



\end{document}